\documentclass{amsart}
\usepackage{amsmath}
\usepackage{amscd}
\usepackage{amssymb}

\theoremstyle{plain}
\newtheorem{thm}{Theorem}
\newtheorem{theorem}[thm]{Theorem}
\newtheorem{lemma}[thm]{Lemma}
\newtheorem{corollary}[thm]{Corollary}
\newtheorem{proposition}[thm]{Proposition}

\theoremstyle{definition}

\newtheorem{remark}[thm]{Remark}

\newtheorem{defn-thm}[thm]{Definition-Theorem}

\newcommand{\End}{{\text{End}}}
\newcommand{\Ad}{{\text{Ad}}}
\newcommand{\bQ}{{\mathbb Q}}
\newcommand{\cF}{{\mathcal F}}
\newcommand{\cN}{{\mathcal N}}
\newcommand{\rc}{{\mathrm c}}
\newcommand{\Lb}{{\mathfrak{b}}}
\newcommand{\Lt}{{\mathfrak{t}}}
\newcommand{\tF}{{\textbf{F}}}
\newcommand{\tk}{{\textbf{k}}}
\newcommand{\Lg}{{\mathfrak g}}
\newcommand{\Lo}{{\mathfrak o}}

\begin{document}
\title[Nilpotent orbits in classical Lie algebras over $\textbf{F}_{2^n}$]{Nilpotent orbits
in classical Lie algebras over $\textbf{F}_{2^n}$ and the Springer
correspondence}
        \author{Ting Xue}
        \address{Department of Mathematics, Massachusetts Institute of Technology,
Cambridge, MA 02139, USA}
        \email{txue@math.mit.edu}

        \begin{abstract}
We give the number of nilpotent orbits in the Lie algebras of
orthogonal groups under the adjoint action of the groups over $\tF_{2^n}$.
Let $G$ be an adjoint algebraic
group of type $B,C$ or $D$ defined over an algebraically closed field of
characteristic 2. We construct the Springer correspondence for the
nilpotent variety in the Lie algebra of $G$.
        \end{abstract}

    \maketitle
\vskip 10pt {\noindent\bf\large Introduction} \vskip 5pt

Throughout this paper, $\tk$ denotes an algebraically closed field
of characteristic 2, $\tF_q$ denotes a finite field of
characteristic 2 and $\bar{\tF}_q$ denotes an algebraic closure of
$\tF_q$.

In \cite{Hes}, Hesselink determines the nilpotent orbits in
classical Lie algebras under the adjoint action of classical groups
over $\tk$. In \cite{Spal}, Spaltenstein gives a parametrization of
these nilpotent orbits by pairs of partitions. We extend Hesselink's
method to study the nilpotent orbits in the Lie algebras of
orthogonal groups over $\tF_q$. Using this extension and
Spaltenstein's parametrization we obtain the number of the nilpotent
orbits over $\tF_q$.

Let $G$ be a connected reductive algebraic group over an
algebraically closed field and $\Lg$ be the Lie algebra of $G$. When
the characteristic of the field is large enough, Springer \cite{Sp1}
constructs representations of the Weyl group of $G$ which are
related to the nilpotent $G$-orbits in $\Lg$. Lusztig \cite{Lu1}
constructs the generalized Springer correspondence which is valid in
all characteristics. Let $G_{ad}$ be an adjoint algebraic group of
type $B,C$ or $D$ over $\tk$ and $\Lg_{ad}$ be the Lie algebra of $G_{ad}$.
 We use a similar construction as in
\cite{Lu1,Lu4} to give the Springer correspondence for $\Lg_{ad}$.
Let $\mathcal{N}_{ad}$ be the set of all pairs
$(\mathrm{c},\mathcal{F})$ where $\mathrm{c}$ is a nilpotent
$G_{ad}$-orbit in $\Lg_{ad}$ and $\mathcal{F}$ is an irreducible
$G_{ad}$-equivariant local system on $\mathrm{c}$ (up to
isomorphism). We construct a bijective map from the set of
isomorphism classes of irreducible representations of the Weyl group
of $G_{ad}$ to the set $\cN_{ad}$. In the case of symplectic group a
Springer correspondence (with a different definition than ours) has
been established in \cite{kato}; in that case centralizers of the
 nilpotent elements are connected \cite{Spal}.
A complicating feature in the orthogonal case is the existence of
non-trivial equivariant local systems on a nilpotent orbit.

\vskip 10pt {\noindent\bf\large Hesselink's classification of
nilpotent orbits over an algebraically closed field} \vskip 5pt

We recall Hesselink's results about orthogonal groups in this
section.

Let $\mathbb{K}$ be a field of characteristic 2, not necessarily
algebraically closed. A form space $V$ is a finite dimensional
vector space over $\mathbb{K}$ equipped with a quadratic form $Q:
V\rightarrow \mathbb{K}$. Let $\langle\cdot,\cdot\rangle: V\times V\rightarrow \mathbb{K}$
be the bilinear form $\langle v,w\rangle=Q(v+w)+Q(v)+Q(w)$. Let
$V^\perp=\{v\in V|\langle v,w\rangle=0,\forall\ w\in V\}$. $V$ is called
non-defective if $V^\perp=\{0\}$, otherwise, it is called defective. $V$
is called non-degenerate if $V\neq\{0\}$, $\text{dim}(V^\perp)\leq
1$ and $Q(v)\neq 0$ for all non-zero $v\in V^{\perp}$.

Let $V$ be a non-degenerate form space of dimension $N$ over
$\mathbb{K}$. Define the orthogonal group $O(V)$ to be
$\{g\in\text{GL}(V)|Q(gv)=Q(v),\ \forall\ v\in V\}$ and define
$\mathfrak{o}(V)$ to be $\{x\in\End (V)|\langle xv,v\rangle=0,\
\forall\ v\in V\text{ and }\text{tr}(x)=0\}$. We write
$O_{N}(\mathbb{K})=O(V)$ and $\Lo_N(\mathbb{K})=\Lo(V)$ when we need
not to specify $V$.  In the case where $\mathbb{K}$ is algebraically
closed let $SO(V)$ be the identity component of $O(V)$ and write
$SO_{N}(\mathbb{K})=SO(V)$.

An element in $\mathfrak{o}(V)$ is nilpotent if and only if it is
nilpotent in $\End (V)$. Let $T$ be a nilpotent element in $\Lo(V)$.
There exists a unique sequence of integers $p_1\geq\cdots\geq
p_s\geq 1$ and a family of vectors $v_1,\ldots,v_s$ such that
$T^{p_i}v_i=0$ and the vectors $T^{q_i}v_i$, $0\leq q_i\leq p_i-1$
form a basis of $V$. We write $p(V,T)=(p_1,\ldots,p_s)$. Define the
index function $\chi(V,T): \mathbb{N}\rightarrow\mathbb{Z}$ by
$\chi(V,T)(m)=\min\{k\geq 0|T^mv=0\Rightarrow Q(T^kv)=0\}.$

Define a form module to be a pair $(V,T)$ where $V$
is a non-degenerate form space and $T$ is a nilpotent element in $
\mathfrak{o}(V)$. To study the nilpotent $O(V)$-orbits in $\mathfrak{o}(V)$ is
equivalent to classify the form modules $(V,T)$ on the form space
$V$. Let $A=\mathbb{K}[[t]]$ and regard $(V,T)$ as an
$A$-module by $(\sum a_nt^n)v=\sum a_nT^nv$. In order to classify
the form modules, Hesselink identifies a form module $(V,T)$ with an
abstract form module $(V,\varphi,\psi)$ (for definition see
\cite{Hes}) via $\varphi: V\times V\rightarrow E,\ (v,w)\mapsto\sum
\langle t^nv,w\rangle t^{-n}$ and $\psi: V\rightarrow E_0,\ v\mapsto\sum
Q(t^nv)t^{-2n}$, where $E$ is the vector space spanned by the linear
functionals $t^{-n}: A\rightarrow\mathbb{K},\ \sum a_it^i\mapsto
a_n,\ n\geq 0$, and $E_0$ is the subspace $\sum_{n\geq 0}\mathbb{K} t^{-2n}$. $E$ is
considered as an $A$-module by $(au)(b)=u(ab)$ for $a,b\in A,u\in
E$. We write $V=(V,\varphi,\psi)$ for simplicity. Define $\mu(v)=\min\{k\geq 0|t^kv=0\}$ for
an element $v$ in the $A$-module $V$ (or $E$).

A form module $V$ is called indecomposable if for every orthogonal
decomposition $V=V_1\oplus V_2$ we have $V_1=0$ or $V_2=0$. Every
form module $V$ has an orthogonal decomposition into indecomposable
submodules $V=\bigoplus_{i=1}^k V_i$. To classify the form modules,
the indecomposable ones are classified first.
\begin{proposition}[\cite{Hes}]\label{prop-indec} Let $V$
be a non-degenerate indecomposable form module. There exist
$v_1,v_2\in V$ such that $V=Av_1\oplus Av_2$ and
$\mu(v_1)\geq\mu(v_2)$. Put
$m=\mu(v_1),m'=\mu(v_2),\Phi=\varphi(v_1,v_2)$ and
$\Psi_i=\psi(v_i)$. One of the following conditions holds:
\begin{enumerate}
\item[(i)] $m'=\mu(\Phi)=m$, $\ \mu(\Psi_i)\leq 2m-1;$
\item[(ii)] $m'=\mu(\Phi)=m-1,\ \mu(\Psi_1)=2m-1>\mu(\Psi_2).$
\end{enumerate}
Conversely, let $m\in\mathbb{N}$, $m'\in\mathbb{N}\cup\{0\}$,
${\Phi\in E}$, $\Psi_1,\Psi_2\in E_{0}$
 be given satisfying $\mathrm{(i)}$ or $\mathrm{(ii)}$.
 Up to a canonical isomorphism there exists a
unique form module $V=Av_1\oplus Av_2$ with
$m=\mu(v_1),m'=\mu(v_2),\Phi=\varphi(v_1,v_2)$ and
$\Psi_i=\psi(v_i)$. This form module is non-degenerate and
indecomposable.
\end{proposition}

Now assume $\mathbb{K}$ is algebraically closed. The form modules
over $\mathbb{K}$ are classified as follows. Let $V=(V,T)$ be a
non-degenerate from module with
$p(V,T)=(\lambda_1,\ldots,\lambda_1,\ldots,\lambda_k,\ldots,\lambda_k)$
where $\lambda_1>\cdots>\lambda_k\geq 1$ and index function
$\chi=\chi(V,T)$. Let $m_i\in\mathbb{N}$ be the multiplicity of
$\lambda_i$ in $p(V,T)$. The isomorphism class of $V$ is determined
by the symbol
$$S(V,T)=(\lambda_1)_{\chi(\lambda_1)}^{m_1}(\lambda_2)_{\chi(\lambda_2)}^{m_2}
\cdots(\lambda_k)_{\chi(\lambda_k)}^{m_k}.$$ A symbol $S$ is the symbol of an isomorphism
class of non-degenerate form modules if and only if it satisfies the
following conditions

(i) $\chi(\lambda_i)\geq\chi(\lambda_{i+1})$ and
$\lambda_i-\chi(\lambda_i)\geq \lambda_{i+1}-\chi(\lambda_{i+1})$,
for $i=1,\ldots,k-1$;

(ii) $\frac{\lambda_i}{2}\leq\chi(\lambda_i)\leq\lambda_i$, for $i=1,\ldots,k$;

(iii) $\chi(\lambda_i)=\lambda_i$ if $m_i$ is odd, for $i=1,\ldots,k$;

(iv) $\{\lambda_i|m_i\text{ odd }\}=\{m,m-1\}\cap\mathbb{N}$
for some $m\in\mathbb{Z}$.

In the following we denote by a symbol either a form module in the
isomorphism class or the corresponding nilpotent orbit.

\vskip 10pt {\noindent\bf\large Isomorphism classes of form modules
and nilpotent orbits over $\tF_q$} \vskip 5pt

Note that the classification of the indecomposable modules
(Proposition \ref{prop-indec}) is still valid. Similar to \cite{Hes}
Section 3.5, we first normalize the non-degenerate indecomposable
form modules over $\tF_q$.
\begin{proposition}\label{prop-nind}
Fix an element $\delta$ in $\tF_q$ such that
$\delta\notin\{x+x^2|x\in\tF_q\}$. The non-degenerate indecomposable
form modules over $\tF_q$ are

$\mathrm{(i)}$ $W_l^0(m)=Av_1\oplus Av_2,\ [\frac{m+1}{2}]\leq l\leq
m$, with
 $\mu(v_1)=\mu(v_2)=m$, $\psi(v_1)=t^{2-2l},\psi(v_2)=0$ and $\varphi(v_1,v_2)=t^{1-m}$;

$\mathrm{(ii)}$ $W_l^\delta(m)=Av_1\oplus Av_2,\ \frac{m+1}{2}\leq
l\leq m$, with $\mu(v_1)=\mu(v_2)=m$, $\psi(v_1)=t^{2-2l}$,
$\psi(v_2)=\delta t^{2l-2m}$ and $\varphi(v_1,v_2)=t^{1-m}$;

$\mathrm{(iii)}$ $D(m)=Av_1\oplus Av_2$ with
$\mu(v_1)=m,\mu(v_2)=m-1$, $\psi(v_1)=t^{2-2m}$, $\psi(v_2)=0$ and
$\varphi(v_1,v_2)=t^{2-m}$.
\end{proposition}

Let $V$ be a non-degenerate form space over $\bar{\tF}_q$. An
isomorphism class of form modules on $V$ over $\bar{\tF}_q$ may
decompose into several isomorphism classes over
 $\tF_q$.
\begin{proposition}\label{prop-2}
Let $W$ be a form module
$(\lambda_1)_{\chi(\lambda_1)}^{m_1}(\lambda_2)_{\chi(\lambda_2)}^{m_2}
\cdots(\lambda_k)_{\chi(\lambda_k)}^{m_k}$ on the form space $V$.
Denote by $n_1$ the cardinality of $\{1\leq i\leq
k-1|\chi(\lambda_i)+\chi(\lambda_{i+1})\leq \lambda_i,\
\chi(\lambda_i)\neq\lambda_i/2\}$ and by $n_2$ the cardinality of
$\{1\leq i\leq k|\chi(\lambda_i)+\chi(\lambda_{i+1})\leq \lambda_i,\
\chi(\lambda_i)\neq\lambda_i/2\}$ (here define $\chi(\lambda_{k+1})=0$).

$\mathrm{(i)}$ If $V$ is defective, the isomorphism class of $W$ over $\bar{\tF_q}$
decomposes into $2^{n_1}$ isomorphism classes over $\tF_q$.

$\mathrm{(ii)}$ If $V$ is non-defective, the isomorphism class of
$W$ over $\bar{\tF_q}$ decomposes into $2^{n_2}$ isomorphism classes over $\tF_q$.
\end{proposition}

Note that we have two types of non-defective form spaces of
dimension $2n$ over $\tF_q$, $V^+$ with a quadratic form of Witt
index $n$ and $V^{-}$ with a quadratic form of Witt index $n-1$. We
denote $O(V^+)$ ($O(V^-)$) by $O^{+}_{2n}(\tF_q)$ (
$O^{-}_{2n}(\tF_q)$) and $\Lo(V^+)$ ($\Lo(V^-)$) by
$\mathfrak{o}^{+}_{2n}(\tF_q)$ ($\mathfrak{o}^{-}_{2n}(\tF_q)$)
respectively. Let $SO^+_{2n}(\tF_q)=O^{+}_{2n}(\tF_q)\cap
SO_{2n}(\bar{\tF}_q)$.

\begin{corollary}\label{coro-n}
$\mathrm{(i)}$ The nilpotent $O_{2n+1}(\bar{\tF}_q)$-orbit
$(\lambda_1)_{\chi(\lambda_1)}^{m_1}\cdots(\lambda_k)_{\chi(\lambda_k)}^{m_k}
$ in $\mathfrak{o}_{2n+1}(\bar{\tF}_q)$ decomposes into $2^{n_1}$
$O_{2n+1}(\tF_q)$-orbits in $\mathfrak{o}_{2n+1}(\tF_q)$.

$\mathrm{(ii)}$ If $\chi(\lambda_i)=\lambda_i/2\ ,i=1,\ldots,k$, the
nilpotent $O_{2n}(\bar{\tF}_q)$-orbit
$(\lambda_1)_{\chi(\lambda_1)}^{m_1}\cdots(\lambda_k)_{\chi(\lambda_k)}^{m_k}$
in $\mathfrak{o}_{2n}(\bar{\tF}_q)$ remains one
$O_{2n}^+({\tF}_q)$-orbit in $\mathfrak{o}^+_{2n}(\tF_q)$;
otherwise, it decomposes into $2^{n_2-1}$ $O^+_{2n}({\tF}_q)$-orbits in
$\mathfrak{o}^+_{2n}(\tF_q)$ and $2^{n_2-1}$
$O^-_{2n}({\tF}_q)$-orbits in $\mathfrak{o}^-_{2n}(\tF_q)$.

Here $n_1,n_2$ are as in Proposition \ref{prop-2}.
\end{corollary}
\begin{remark}
If $\chi(\lambda_i)=\lambda_i/2\ ,i=1,\ldots,k$, then $n$ is even.
If $\chi(\lambda_i)\neq\lambda_i/2$ for some $i$, then $n_2\geq 1$.
\end{remark}

Using Corollary \ref{coro-n},
we can give a bijective proof of the following proposition.

\begin{proposition}\label{prop-num}
$\mathrm{(i)}$ The number of nilpotent $O_{2n+1}(\tF_q)$-orbits in
$\mathfrak{o}_{2n+1}(\tF_q)$
 is $p_2(n)$.

$\mathrm{(ii)}$ The number of nilpotent $O_{2n}^+(\tF_q)$-orbits in
$\mathfrak{o}_{2n}^+(\tF_q)$
 is $\frac{1}{2}p_2(n)$ if $n$ is odd and is
$\frac{1}{2}(p_2(n)+p(\frac{n}{2}))$ if $n$ is even.

Here $p_2(k)$ is the number of pairs of partitions $(\alpha,\beta)$
such that $|\alpha|+|\beta|=k$ and $p(k)$ is the number of
partitions of the integer $k$.
\end{proposition}

\begin{corollary}\label{coro}
The number of nilpotent $SO_{2n}^+(\tF_q)$-orbits in
$\mathfrak{o}_{2n}^+(\tF_q)$ is $\frac{1}{2}p_2(n)$ if $n$ is odd
and is $\frac{1}{2}p_2(n)+\frac{3}{2}p(\frac{n}{2})$ if $n$ is even.
\end{corollary}

Let $G=SO_N(\tk)$. The Lie algebra $\Lg$ of $G$ is $\Lo_N(\tk)$. Let
$\mathcal{N}$ be the set of all pairs $(\mathrm{c},\mathcal{F})$
where $\mathrm{c}$ is a nilpotent
$G$-orbit in $\Lg$ and $\mathcal{F}$ is an irreducible
$G$-equivariant local system on $\mathrm{c}$ (up to
isomorphism). Let $a$ be the number of irreducible representations of the Weyl
 group $W$ of $G$. We show that the number of elements in $\cN$ is equal to $a$.
To see this we can assume $\tk=\bar{\tF}_2$. In this case for $q$ a power of $2$,
let $G(\tF_q)$, $\mathfrak{g}({\tF}_q)$ be the fixed points of a split Frobenius
map $\mathfrak{F}_q$ relative to $\tF_q$ on $G$, $\Lg$. From Proposition
\ref{prop-num} (i) and Corollary \ref{coro} we see that the number of
nilpotent $G(\tF_q)$-orbits in
$\mathfrak{g}({\tF}_q)$ is equal to $a$.
Pick representatives $x_1,\cdots,x_M$ for
the nilpotent $G$-orbits in $\Lg$. If $q$ is large enough, the
Frobenius map $\mathfrak{F}_q$ keeps $x_i$ fixed
and acts trivially on $Z_{G}(x_i)/Z_{G}^0(x_i)$. Then the number of
$G(\tF_q)$-orbits in the $G$-orbit of $x_i$ is equal to the number
of irreducible representations of $Z_{G}(x_i)/Z_{G}^0(x_i)$ hence to
the number of $G$-equivariant irreducible local systems on the
$G$-orbit of $x_i$.

Assume $G_{ad}$ is an adjoint
group over $\tk$ of the same type as $G$ and $\Lg_{ad}$ is the Lie
algebra of $G_{ad}$. Let $G_{ad}(\tF_q)$, $\Lg_{ad}(\tF_q)$ be defined like
$G(\tF_q)$, $\Lg(\tF_q)$. The Lie algebra $\Lg_{ad}$
is not isomorphic to $\Lg$. But the number of
nilpotent $G_{ad}(\tF_q)$-orbits in $\Lg_{ad}(\tF_q)$
 is the same as the number of nilpotent $G(\tF_q)$-orbits in
$\Lg(\tF_q)$. In fact, we have a morphism $G\rightarrow G_{ad}$
which is an isomorphism of abstract groups and an obvious bijective
morphism $\mathcal{U}\rightarrow \mathcal{U}_{ad}$ between the
nilpotent variety $\mathcal{U}$ of $ \Lg$ and the nilpotent variety
$\mathcal{U}_{ad}$ of $\Lg_{ad}$. Thus the nilpotent orbits in
$\Lg$ and $\Lg_{ad}$ are in bijection and the corresponding
 component groups of centralizers are isomorphic. It follows that the
number of elements in $\cN_{ad}$ (as in the introduction)
is equal to the number of elements
in $\cN$.

Note that the argument in the last two paragraphs also applies for
the symplectic group.

\vskip 10pt {\noindent\bf\large Springer correspondence} \vskip 5pt

Assume $G_{ad}$ is an adjoint group of type $B_r$, $C_r$ or $D_r$ over
$\tk$ and $\Lg_{ad}$ is the Lie algebra of $G_{ad}$. Let $\cN_{ad}$ be as
in the introduction. We give the
Springer correspondence for $\Lg_{ad}$. The following lemma plays an
important role in our construction.
\begin{lemma}\label{l-1}
There exist regular semisimple elements in $\Lg_{ad}$ and they form
an open dense subset in $\Lg_{ad}$.
\end{lemma}
\begin{remark}
This lemma is not always true if the group is not adjoint.
\end{remark}

Fix a Borel subgroup $B$ and a maximal torus $T\subset B$ in
$G_{ad}$. Let $W=N_{G_{ad}}(T)/T$ be a Weyl group of $G_{ad}$.
Denote the Lie algebra of $B$ by $\Lb$ and the Lie algebra of $T$ by
$\Lt$. Let $\Lt_0$ be the set of regular elements in $\Lt$ and $Y$
be the set of regular semisimple elements in $\Lg_{ad}$. By Lemma
\ref{l-1}, the closure $\bar{Y}$ of $Y$ in $\Lg_{ad}$ is $\Lg_{ad}$.
Let $\widetilde{Y}=\{(x,gT)\in Y\times G_{ad}/T|\Ad(g^{-1})(x)\in
\Lt_0\}$. Define $\pi: \widetilde{Y}\rightarrow Y$ by $\pi(x,gT)=x$.
We have that $\pi$ is a principal $W$-bundle, hence
$\pi_!\bar{\mathbb{Q}}_{l\widetilde{Y}}$ is a well defined local
system on $Y$ and thus the intersection cohomology complex
$IC(\Lg_{ad},\pi_!\bar{\mathbb{Q}}_{l\widetilde{Y}})$ is well
defined.
 Let $X=\{(x,gB)\in\Lg_{ad}\times G_{ad}/B|\Ad(g^{-1})x\in\Lb\}$. Define
$\varphi:X\rightarrow\Lg_{ad}$ by $\varphi(x,gB)=x$.
\begin{proposition}\label{p-2}
$\varphi_!\bar{\mathbb{Q}}_{lX}$ is canonically isomorphic to
$IC(\Lg_{ad},\pi_!\bar{\mathbb{Q}}_{l\widetilde{Y}})$.
\end{proposition}
We have
$\End(\varphi_!\bar{\mathbb{Q}}_{lX})=
\End(\pi_!\bar{\mathbb{Q}}_{l\widetilde{Y}})=\bar{\bQ}_l[W]$.
Let $\hat{W}$ be a set of representatives for the isomorphism
classes of simple $W$-modules. We have a canonical decomposition $
\varphi_!\bar{\bQ}_{lX}=\bigoplus_{\rho\in
\Hat{W}}(\rho\otimes(\varphi_!\bar{\bQ}_{lX})_{\rho}) $.
Set $\bar{Y}^\omega=\{x\in\bar{Y}|x\text{ nilpotent}\}$.

\begin{theorem}\label{mp-1}
Let $d_0=\dim G_{ad}-\dim T$. For any $\rho\in \hat{W}$, there is a
unique $(\mathrm{c},\mathcal{F})\in\mathcal{N}_{ad}$ such that
$(\varphi_!\bar{\bQ}_{lX})_\rho|_{\bar{Y}^\omega}[d_0]$ is
$IC(\bar{\rc},\cF)[\dim\rc]$ regarded as a simple perverse sheaf on
$\bar{Y}^\omega$ (zero outside $\bar{\rc}$), where $\bar{\rc}$ is the closure
of $\rc$ in $\bar{Y}^\omega$. Moreover, $\rho\mapsto
(\rc,\cF)$ defines a bijective map $\gamma:
\hat{W}\rightarrow\mathcal{N}_{ad}$.
\end{theorem}

A corollary is that in this case there are no cuspidal local systems
similarly defined as in \cite{Lu1}. This result does not extend to exceptional
Lie algebras. (In type $F_4$, characteristic 2, the results of \cite{Spa2} suggest
that a cuspidal local system exists on a nilpotent class.)

\vskip 10pt {\noindent\bf\large Acknowledgement} \vskip 5pt

I would like to thank Professor George Lusztig for his guidance and
encouragement during this research.

\end{document}